%% file: hyperfiber17.tex
\documentclass[12pt,reqno]{amsart}
\usepackage{geometry}       
\usepackage{graphicx}
\usepackage{amssymb}
\usepackage{epstopdf}

\usepackage{epsfig}

\usepackage{amsthm}

\def\mathscr{\mathcal}
\usepackage[all]{xy}

\theoremstyle{plain}
\newtheorem{theorem}[equation]{Theorem}
\newtheorem{proposition}[equation]{Proposition}
\newtheorem{corollary}[equation]{Corollary}
\newtheorem{lemma}[equation]{Lemma}

\theoremstyle{definition}
\newtheorem{definition}[equation]{Definition}
\newtheorem{remark}[equation]{Remark}
\newtheorem{example}[equation]{Example}
\numberwithin{equation}{section}

\newcommand{\Z}{{\mathbb Z}}

\renewcommand{\leq}{\leqslant}
\renewcommand{\geq}{\geqslant}

\usepackage{color}
\begin{document}

\title[Infinite Fibering]{Hyperbolic Groups Which Fiber in Infinitely Many Ways}

\author[T. Mecham, A. Mukherjee]{TaraLee Mecham, Antara Mukherjee}
\address{Dept.\ of Mathematics\\
        University of Oklahoma\\
        Norman, OK 73019\\
        USA}
\email{tmecham@math.ou.edu, amukherjee@ou.edu}
\date{}

\begin{abstract}
We construct examples of free-by-cyclic hyperbolic groups which
fiber in infinitely many ways over $\Z$.  The construction involves adding a
specialized square 2-cell to a non-positively curved, squared
2-complex defined by labeled oriented graphs. The fundamental groups of the resulting complexes are hyperbolic, free-by-cyclic and can be mapped onto $\Z$ in infinitely many ways.
\end{abstract}

\maketitle

\section{Introduction}
Topologists have been interested in fiberings of 3-manifolds since
the early 1960's. A 3-manifold is said to fiber over $S^1$ if it
is the mapping torus of a surface via a homeomorphism. A key to
relating this topological fibering to the notion of algebraic
fibering was expressed by Stallings in 1962 \cite{St}. He
proved:
   \emph{For an irreducible compact 3-manifold, $M^3$, if
   $\pi_1(M^3)$ has a finitely generated, normal subgroup G (which is not $\Z_2$), whose quotient group is
   $\Z$, then $M^3$ fibers over $S^1$. The fiber here is a 2-manifold $T$ embedded in $M^3$, such that $\pi_1(T)= G$. }
    The fundamental groups of these manifolds are a special class of HNN-extensions that extend the fundamental group of the surface by $\Z$.

     In 1969, Tollefson showed that there were 3-manifolds that fibered over $S^1$ in
infinitely many ways, \cite{To}. From the group theory
perspective, this means that the fundamental groups of  these
3-manifolds can be expressed in infinitely many ways as
surface-by-cyclic groups with the genus of the surface going to
infinity. In 1998 Thurston \cite{Wi}
 proved that \emph{every atoroidal three-manifold that
fibers over the circle has a hyperbolic structure}.  Using branched covers of the 3-torus, Hilden, Lozano
and Montesinos-Amilibia \cite{HiLo}
     constructed the first explicit examples of hyperbolic manifolds which have infinitely many fibrations over
$S^1$, in 2006.
 The main theorem of this paper produces
analogous results in the context of hyperbolic, free-by-cyclic
groups.

\begin{theorem}
\label{main} There exist hyperbolic, free-by-cyclic groups which
have  infinitely many epimorphisms onto $\Z$. Furthermore, each of
these groups is isomorphic to $F_n \rtimes \Z$ for infinitely many
$n$.
 \end{theorem}
 The groups we produce will be the first concrete examples of hyperbolic, free-by-cyclic groups which fiber in infinitely many ways over $\Z$.  We prove Theorem \ref{main} by constructing groups
that are hyperbolic, free-by-cyclic and map onto $\Z^2$. These
groups are the fundamental groups of specially constructed
non-positively curved, squared 2-complexes.

 The paper is organized as follows:
 Section 2 gives some
background about non-positively curved squared 2-complexes,
$\delta$-hyperbolic spaces, and Morse Theory  necessary for the
construction of the groups. Section 3 gives a description of the
building blocks for the construction in Section 4. We
give the actual construction of groups with reducible monodromy automorphism,  proving Theorem \ref{main},  in Section 4.
Section 5 provides an alternative construction to produce an
example of a group which may have an irreducible automorphism.
This construction does not use the building blocks from Section 3.  We conclude the paper with
Section 6, where we present some questions that are yet to be
answered.

This paper would not have been possible without the support and guidance of our advisor, Dr. Noel Brady.


\section{Background and Definitions}

Theorem \ref{main} requires that we produce free-by-cyclic, hyperbolic groups, so we will recall some facts about
$\delta$-hyperbolic spaces and combinatorial Morse Theory in this section.
 Since the  groups we construct are the fundamental groups of
non-positively curved, squared 2-complexes, we first recall some
facts about non-positively curved spaces.  Then we quote Gromov's
Flat Plane Theorem in \cite{Gr} which helps us show that
the groups are indeed hyperbolic.
 Finally, we will present
some facts about Morse Theory and conclude the section with
Proposition \ref{FbyC} which we will use
to both prove that our groups are free-by-cyclic, and compute the free rank of the groups.

We start with some definitions associated with non-positively curved spaces.
These definitions can be found in Bridson and
Haefliger's  book \cite{BrHa}.

  \begin{definition}[Piecewise  2-complex]
  \label{PE2cmplx}
    A \emph{piecewise Euclidean 2-complex} $K$ is
    a 2-complex obtained from a collection of convex cells in the Euclidean plane
    by identifying their edges via isometries.
\end{definition}

  \begin{definition}[Link of a vertex]
\label{link}
 The link of a vertex $v$ in a convex 2-cell $c$
(denoted by $Lk(v,c)$) of a 2-complex, is defined as the
collection of unit tangent vectors to the space at $v$ which point
to $c$. This is a spherical cell of dimension $(n-1)$.
  \end{definition}

Suppose there are two convex 2-cells $c, c'$ with vertices $v, v'$
respectively which are identified along the faces via an isometry
$f$ which identifies $v$ with $v'$.  The derivative of $f$ gives a
spherical isometry which identifies a face of $Lk(v,c)$ and a face
of $Lk(v',c')$. Thus, the the link of a vertex in a 2-complex is
naturally described as a piecewise spherical complex.

  \begin{definition}[Link Condition]
  \label{Linkcon}
    The link of a vertex satisfies the
     \emph{large link condition} if it has no circuits of length $< 2\pi $.  This definition is the same as saying that  the link
     satisfies the no empty triangles condition or equivalently,  the link is a flag
     complex.   See \cite{BrHa} for more details.

  \end{definition}

  \begin{definition}[Non-positively curved 2-complex]
  \label{NPC2cmplx}
     A piecewise Euclidean 2-complex $K$ is said to be \emph{non-positively curved} if $Lk(v,K)$ is large (i.e., no circuits of length less than $2\pi$) for all $v \in K^{(0)}$.
      In this context a link is a simplicial complex, which implies it has no edge loops and no bigons.
  \end{definition}

  The following is a theorem stated by Gromov  \cite{Gr}, which was later proved by Bridson
\cite{Br}.  We will use this theorem to show that the universal
covers of the 2-complexes we construct are hyperbolic, which will
imply that the corresponding fundamental groups are hyperbolic. This
result is analogous to  Eberlein's \cite{Eb} earlier work, where
he showed that:  \emph{A compact non-positively curved Riemannian
manifold either has a word-hyperbolic fundamental group or its
universal cover contains an isometrically embedded copy of
the Euclidean plane}.

  \begin{theorem}[The Flat Plane Theorem \cite{BrHa}]
  \label{Hyp}
     A proper cocompact CAT(0) space is \emph{hyperbolic} if and only if it does not contain a subspace
     isometric to $\mathbb E^2$.
  \end{theorem}

 We will also need to
ensure that our groups are free-by-cyclic.  We will use
combinatorial Morse Theory to establish this fact.  Following are
the basic definitions necessary for this method of proof.  These
definitions, and proofs for the theorems, can be found in \cite{BaBr}.

  \begin{definition}[Affine cell-complex]
  \label{Affine}
    A finite-dimensional cell-complex $X$ is said to be an \emph{affine cell-complex} if it has the following
    structure.  An integer $ m > dim(X)$ is given, and for each cell $e$ of $X$ we are given a convex
    polyhedral cell $C_e \subset$ $\mathbb R^m$ and a characteristic function $\chi_e:C_e\rightarrow e$
    such that the restriction of $\chi_e$ to any face of $C_e$ is a characteristic function of another cell,
    possibly pre-composed by a partial affine homeomorphism of $\mathbb R^m$.
  \end{definition}

  In our construction, our highest dimensional cell will be two-dimensional.  In fact, the complexes will be  made of Euclidean squares identified
  along the edges, so they will be affine.

     \begin{definition}[Morse Function]
  \label{Morse}
    A map $f:X \rightarrow \mathbb R$ defined on an affine cell-complex $X$ is a \emph{Morse
function}
    if:

        \begin{itemize}
             \item
                 for every cell $e$ of $X$, $f_{\chi_e} :C_e \rightarrow \mathbb R$ extends to an
affine map $
                 \mathbb R^m \rightarrow \mathbb R$ and $f_{\chi_e}$ is constant only when $dim(e)
= 0$, and

            \item
                the image of the $0$-skeleton is discrete in $\mathbb R$.

        \end{itemize}
  \end{definition}

  One can think of a Morse function as a height function where each cell of the complex has  a unique maximum and a
  unique minimum.  Since we want to show that our groups are free-by-cyclic, we will be interested in maps
  from our complexes to the circle.  The following definition relates these maps to Morse functions.

   \begin{definition}[Circle-valued Morse function]
   \label{CircleMorse}
      A \emph{circle-valued Morse function} on an affine cell complex $X$ is a cellular map $f:X
\rightarrow
      S^1$, with the property that $f$ lifts to a Morse function between universal covers.  When the context is clear, we will simply refer to these as Morse functions.

   \end{definition}

  The links of the vertices of a 2-complex are graphs.
  It turns out that in the presence of a Morse function,  certain subgraphs of the link give us important information about the structure of the
   fundamental group of the 2-complex.  The following definition and proposition describe this relationship.

  \begin{definition}[Ascending and descending links]
  \label{AscDesLinks}
     Suppose $X$ is an affine cell complex and $f:X \rightarrow S^1$ is a circle-valued Morse
function.
     Choose an orientation of $S^1$,  which lifts to one of $\mathbb R$, and lift $f$ to a map of
universal
     covers $\widetilde{f}:\widetilde{X}\rightarrow \mathbb R$.  Let $v \in X^{(0)} $, and note
that the link of
      $v $ in $X$ is naturally isomorphic to the link of any lift $\widetilde{v}$ of $v$ in
$\widetilde{X}$.  We
      say that a cell $\widetilde{e}\subset \widetilde{X}$ contributes to the \emph{ascending}
(respectively
      \emph{descending}) link of $\widetilde{v}$ if $ \widetilde{v} \in \widetilde{e}$ and only if
$\widetilde{f}|
      _{\widetilde{e}}$ achieves its minimum (respectively maximum) value at $\widetilde{v}$.  The
      ascending (respectively descending) link of $v$ is then defined to be the subset of
$Lk(v,X)$
      naturally identified with the ascending (respectively descending) link of $\widetilde{v}$.

  \end{definition}

The following proposition is an important result of Morse Theory and
helps us check that the fundamental groups of the 2-complexes we
construct are indeed free-by-cyclic and gives us a formula for determining the free rank of the groups.
  \begin{proposition}[Free-by-Cyclic \cite{BaBr}]
  \label{FbyC}
     If $f:X \rightarrow S^1$ is a circle-valued Morse function on the 2-complex X all of whose
ascending
     and descending links are trees, then X is aspherical, and $\pi_1(X)$ is free-by-cyclic.  In
fact, $\pi_1
     (X) = F_k \rtimes \Z$ where $k = 1 - \chi (f^{-1}(p))$ for $p\in S^1$.
 \end{proposition}


 \section{Building Blocks}
\label{blocks}

The groups that we construct
 in Section 4 are one-relator quotients of free
products of some building block groups.  In this section we describe
these building block groups and show that they are the fundamental
groups of non-positively curved 2-complexes.  We also show that these
groups are hyperbolic and free-by-cyclic.

\begin{remark}  We will be looking at 2-complexes which have a
unique vertex for the rest of the paper.
\end{remark}

We start this section by defining a Labeled Oriented Graph (LOG).
The group presentation corresponding to a particular LOG is used as the building block group for the example
constructed in section 4.  The following definition is taken from
\cite{BaBr}.

  \begin{definition}[Labeled Oriented Graph]
  \label{LOG}
    A \emph{labeled oriented graph} or LOG, consists of a finite, directed graph with distinct
labels on all
    the vertices and oriented edge labels taken from the set of vertex labels.  An LOG defines a
finite
    group presentation in the following manner.  The set of generators is in 1-1 correspondence
with the
    set of vertex labels and the set of relators is in 1-1 correspondence with the set of edges,
so that an
    edge labeled \emph{a} oriented from vertex \emph{u} to vertex \emph{v} corresponds to a
conjugation
    relation $ava^{-1}=u$.  In the case that the graph is a tree,  we call it a \emph{labeled
oriented tree} or
    LOT, and call the corresponding finite presentation an LOT presentation.
  \end{definition}

 Howie \cite{Ho} was one of the first people to define and study the group
presentation of LOTs.


\begin{figure}[htbp] \begin{center} 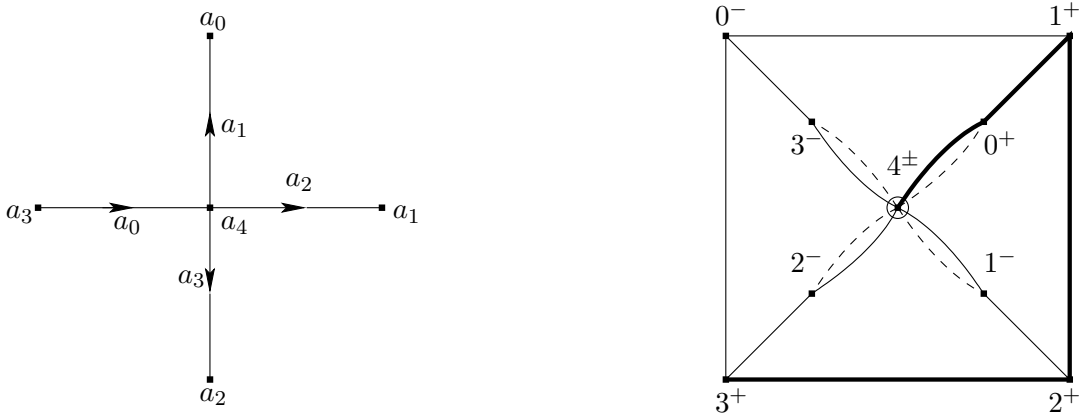 \caption{The LOT presentation of $L_{(a_i)}$  and the link of $v$ in $X_{L_{(a_i)}}$.  The circle in the center of the link represents the points $4^+$ and $4^-$.  The dashed lines are edges connecting to $4^-$.  There is no edge connecting $4^+$ to $4^-$.    The descending link is shown in bold.} \label{lotandlink} \end{center} \end{figure}


  \begin{remark}
    The main component of our construction is an LOT group with five generators
$a_0,...,a_4$
    denoted by $L_{(a_i)}$.   Figure \ref{lotandlink} shows the graph of the LOT.
The corresponding presentation 2-complex  $X_{L_{(a_i)}}$, consists of  a bouquet of five
circles and four square 2-cells. The group presentation of this complex is given by,
  $$
    L_{(a_i)} = \langle a_0,a_1,a_2,a_3,a_4 \hspace{0.2cm}|\hspace{0.2cm} a_4 = a_0^{a_1},\hspace
    {0.2cm} a_4 = a_1 ^{a_2},\hspace{0.2cm} a_4 = a_2 ^{a_3},\hspace{0.2cm}a_3 = a_4^{a_0}\rangle
  $$
where $a_i^{a_j}= a_ja_ia_j^{-1}$.

The link of the vertex of the corresponding presentation 2-complex is also shown in Figure \ref{lotandlink}.  $X_{L_{(a_i)}}$ is a piecewise Euclidean 2-complex and it is easy to
see that the link is large and hence, the universal cover
$\widetilde{X}_{L_{(a_i)}}$ is CAT(0). In order to show that
$L_{(a_i)}$ is hyperbolic we will need to use the Flat Plane
Theorem.   The following definition will be useful in its
application.

\end{remark}

\begin{definition}[Poison Corners]
\label{poison} A \emph{poison corner} is a corner of a 2-cell for
which the corresponding edge in the link is not a part of any
length-four circuit.
\end{definition}

\begin{remark}
This definition is helpful to us because the 2-cells of the
presentation 2-complex $X_{L_{(a_i)}}$ are Euclidean square 2-cells.
This means that all edges in the link are of length $\pi/2$. Hence,
any square of our complex that contains a poison corner cannot be
used to develop a flat plane in the universal cover.
\end{remark}

  \begin{theorem}
  \label{XLAI}
    The universal cover of  $X_{L_{(a_i)}}$ has no isometrically embedded flat planes.
  \end{theorem}

\begin{proof}

    The link has no bigons
or empty
   triangles, so the link is large.  We now need only ensure that no planes develop in the universal cover.
    There are four edges of the link that act as bridges from the inner
suspension to
   the outer circuit of the link.  They are:
  $$
       \hspace{0.5cm}a_1^{+} a_0^{+},\hspace{0.5cm}a_2^{+} a_1^{-},\hspace{0.5cm} a_3^{+} a_2^{-},
       \hspace{0.5cm}a_0^{-} a_3^{-}
  $$
 If we remove these four edges, the link reduces to two separate components.  The outer component is the length-four circuit ${a_0}^-
{a_1}^+{a_2}^+{a_3}
   ^+{a_0}^-$.  The combinatorial distance between any two vertices of this component is at least
1.  On
   the other hand, the inner component is the suspension on the two vertices ${a_4}^+$ and
${a_4}^- $
   that connecting them to ${a_0}^+,{a_1}^-,{a_2}^-,$ and ${a_3}^-$.
   The distance between any pair of these four vertices in the inner component is at least 2.  Therefore, any
circuit
   containing one of these bridges has to be of length at least 5.  Hence, the bridges represent poison corners.

   Any isometrically embedded flat plane is tiled by a collection of 2-cells of the presentation
2-complex
   $X_{L_{(a_i)}}$. However, every 2-cell of this complex contains a poison corner.  This means that there are no
isometrically
   embedded flat planes in the universal cover of the complex.
   \end{proof}

  \begin{corollary}
  \label{LAIisHYP}
     $L_{(a_i)}$ is a hyperbolic group.
   \end{corollary}
\begin{proof}
 The universal cover of the 2-complex $X_{L_{(a_i)}}$ is CAT(0) since the link is large and Theorem
\ref{XLAI} further
   implies that the universal cover $\widetilde{X}_{L_{(a_i)}}$ has no isometrically embedded flat
planes, so using the Flat Plane Theorem, we know that $\widetilde{X}_{L_{(a_i)}}$ is $\delta$-hyperbolic and the corresponding group $L_{(a_i)}$ is hyperbolic.
  \end{proof}

Now we will show that $L_{(a_i)}$ is also free-by-cyclic. The
proof of the following proposition is along the same lines as the first part of the proof of Theorem 3 in \cite{BaBr}.

\begin{proposition}
$L_{(a_i)}$ is free-by-cyclic and is isomorphic to $F_4\rtimes \Z$.
\end{proposition}
\begin{proof}
  The presentation 2-complex $X_{L_{(a_i)}}$ of
the group $L_{(a_i)}$, admits a circle-valued Morse function which
is defined by mapping the vertex to the base vertex of the circle
and each directed edge is mapped once around the target circle. The
same map extends linearly over the 2-cells of
$X_{L_{(a_i)}}$.
 The descending link of this Morse function is shown in bold on the link in Figure
 \ref{lotandlink}.
Since both the ascending and descending links are trees, Proposition
\ref{FbyC} implies that the LOT group is free-by-cyclic. The
pre-image of the vertex in the target circle under the Morse
function is a bouquet of four circles, so
the group is isomorphic to  $F_4\rtimes \Z$.
\end{proof}

Although $L_{(a_i)}$ is hyperbolic and free-by-cyclic, it only
fibers in one way over $\Z$ because the relations of the group are
all conjugation relations.  This means that the image of all the
generators will be the same under any map to $\Z$. Therefore there
is only one nontrivial map to $\Z$, up to post composition with a
monomorphism of $\Z$.

\begin{remark}
\label{otherLOT}

It is easy to see that for any LOT group represented by a tree
with $k+1$ vertices, a single valence $k$ vertex, and with edge
and vertex labels permuted as in the $L_{(a_i)}$ case, the
previous theorems hold. In fact, the link will consist of an
outer $k$-gon with $k$ bridges to $k$ inner vertices and the
suspension of these inner vertices on the vertices $a_k^+$ and
$a_k^-$.  This means that the group will be hyperbolic and
isomorphic to $F_k\rtimes \Z$. For arbitrary $k \geq 4$, we will
denote the group and the presentation 2-complex  associated to it
by $L_{(k, a_i)}$, and $X_{L_{(k, a_i)}}$ respectively.

\end{remark}


\section{The Construction of Reducible Examples}
\label{const}

Now we will use the LOT groups from the previous section to prove
Theorem \ref{main}.  We will begin by focusing on the topological
aspects of our construction and show that the groups produced
are free-by-cyclic and map onto $\Z^2$. This will be the content
of Theorem \ref{maps} below.
 Next we will address the geometry of our particular examples by
showing that they are indeed hyperbolic, thus providing the proof
of Theorem \ref{main}.  We will then give a method for computing
the free rank of the kernel of a map from our groups onto $\Z$ in
Lemma \ref{growth}. We conclude this section by writing down
an explicit (reducible) automorphism for the group produced in the
proof of Theorem \ref{main}.
 .

  \begin{remark}
  \label{Notation}
Since all of the spaces that we deal with have a single vertex, we shall denote the fundamental
group
     element corresponding to any 1-cell say $x$, of a 2-complex, by
$[x]$.
If $f$ is the map from a presentation 2-complex, $K$, to a circle
and $f_*$ the map between the corresponding fundamental groups,
i.e, $f_*:\pi _1(K)\rightarrow \Z$, then using the notation given above,
we will define $\alpha = f_*([a])$ and $\beta = f_*([b])$.  We will use this notation throughout the rest of the paper.

  \end{remark}

   \begin{theorem}
   \label{maps}
     Given 2-complexes $K_j$ for $j= 1,2$ satisfying:

        \label[(i)   $|K_j^{(0)}|= 1$,

        \label[(ii)  there exist Morse functions $ f_j: K_j \rightarrow S^1$ such that the
ascending and
        descending links are trees,

        \label[(iii)  there are 1-cells $a_0,a_1 \in K_1$ and $b_0, b_1 \in K_2$ such that
$|\alpha_0 |+
        |\beta_1| = |\alpha_1| + |\beta_0|$,

     then the 2-complex, $K$, obtained from ${K_1 \sqcup K_2}/(K_1^{(0)} \sim K_2^{(0)})$ by
attaching a
     square 2-cell according to the relation $r=a_0^{\pm 1}b_1^{\pm 1}a_1^{\mp1}b_0^{\mp1}$, satisfies:

        \label[(i)   $\pi _1(K)$ is free-by-cyclic and,

        \label[(ii)   $\pi _1(K)$ maps onto $\Z^2$.

   \end{theorem}

\begin{proof}
   \label[(i)
      We start by labeling the unique vertex of $K_j$ by $v_j$ for $j = 1, 2$ and consider the
2-complex  $K_1
      \vee K_2 = {K_1 \sqcup K_2}/v_1 \sim v_2$, calling the common vertex $v$.  The fundamental
group
      of  this complex is the free product of the fundamental groups of $K_1$ and $K_2$.
     By the functorality of the wedge product, $f_1$ and $f_2$ induce the following map:
         $$ f_1\vee f_2 : K_1 \vee K_2 \rightarrow S^1\vee S^1$$
     Using $f_1\vee f_2$, we define a circle-valued Morse function, $pr\circ (f_1\vee f_2)$, on the 2-complex $K_1\vee K_2$ so that $pr\circ (f_1\vee f_2)\circ incl. = f_j$ for $j=1,2$, where $pr$ is the projection map of $S^1\vee S^1$ onto $S^1$. Following is a commutative diagram showing the relationship between these maps.
           \begin{displaymath}
             \xymatrix{K_j \ar[d]_{f_j}\ar@{^{(}->}[r]^{incl.}  &  K_1\vee K_2 \ar[d]^{f_1\vee
f_2} \\ S^1        &       S^1
             \vee S^1\ar[l]^{pr}}
        \end{displaymath}

  The ascending and descending links of $v$ are each 2-component
forests
  (the union of the ascending and descending links of $v_1\in K_1$ and $v_2 \in K_2$).  To
complete the
  construction, we attach a single 2-cell which represents a relation $ r $ of the form $a_0^{\pm 1}b_1^{\pm 1}=
b_0^{\pm 1}a_1^{\pm 1}$ where these are the 1-cells described in
hypothesis  \label[(iii).  The powers in the relation are positive
or negative depending upon the sign of the associated $\alpha _ i$,
$\beta_i$.   Since $\alpha_0 + \beta_1 =  \alpha_1 + \beta_0$, we
can use Lemma 1.5 of Scott and Wall's paper \cite{ScTe} to extend
the Morse function of $K_1\vee K_2$
  over $K$.  We will call this new function $g$.  The extra 2-cell contributes an edge each to the ascending and descending
  links connecting up the two components of the forest in each case.  Thus the ascending and
  descending links are now trees and $\pi_1(K) $ is free-by-cyclic by Proposition \ref{FbyC}.  The
  following commutative diagram shows the Morse function on the complex $K$.

     \begin{displaymath}
          \xymatrix{K_1 \vee  K_2 \ar[d]_{f_1\vee f_2} \ar@{^{(}->}[r]  & K \ar@{.>}[d]^{g}\\
     S^1\vee S^1 \ar[r] _{pr}                             & S^1 }
    \end{displaymath}

\label[(ii)There is also a map from $K_1\vee K_2$ to $S^1\times S^1$ which is shown in the next diagram.  SInce $K_1\vee K_2$ is a
  sub-complex of $K$ we can use the same result as above in \cite{ScTe} and extend this map to a map $g'$ over $K$.

       \begin{displaymath}
        \xymatrix{K_1 \vee  K_2 \ar[d]_{f_1\vee f_2} \ar@{^{(}->}[r]  & K \ar@{.>}[d]^{g'}\\
         S^1\vee S^1  \ar@{^{(}->}[r]   & S^1\times S^1 }
     \end{displaymath}

  This gives us a map  from $\pi_1(K)= \pi_1(K_1) \ast \pi_1(K_2)/<<r>>$ to $\Z^2$.  The map to $\Z^2$ is
  onto since it is composed of $f_{1_*}: \pi_1(K_1) \rightarrow \Z$ and $f_{2_*}: \pi_1(K_2) \rightarrow \Z
  $, where the maps  $f_{j_*}$  are defined by sending generators of $\pi_1(K_j)$ to generators of each
  factor of $\Z$ .

  \end{proof}

  Now all the pieces are in place for us to provide the proof of the main theorem.

\begin{figure}[htbp] \begin{center} 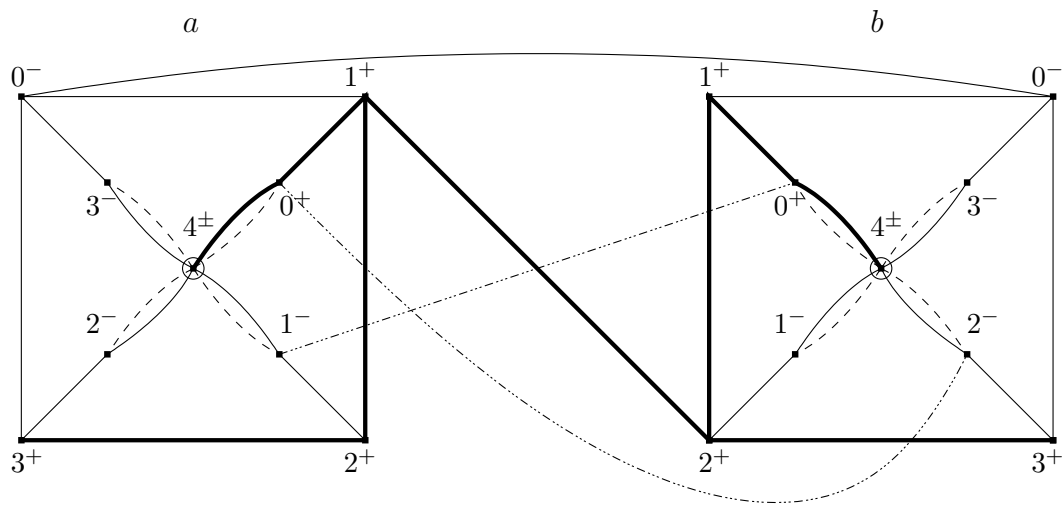 \caption{$Link(v,X_1)$.  The descending link is shown in bold.} \label{lotlotlink} \end{center} \end{figure}

   \begin{proof}[Proof of Theorem \ref{main}]
  \label{LOT-LOT}

    We begin with 2 copies of our LOT group namely, $L_{(a_i)}$, $L_{(b_i)}$ and their associated
    presentation 2-complexes, $X_{L_{(a_i)}}$, $X_{L_{(b_i)}}$.  Using the above construction, we
obtain a new complex $X_1$ such that its fundamental group is:
       $$G_1 = L_{(a_i)}\ast L_{(b_i)}/\langle \langle a_0b_2a_1^{-1}b_0^{-1}\rangle \rangle  =
\langle L_
       {(a_i)}, L_{(b_i)}\hspace{0.2cm} | \hspace{0.5cm}a_0b_2a_1^{-1}b_0^{-1} \rangle $$
  Any choice of relation $a_ib_j = b_ka_l$ would fit the conditions of Theorem \ref{maps}, but not
all of
  these would ensure that our resultant group is hyperbolic, so we have carefully chosen
indices for the generators in the relation so that the link of the
new complex $X_1$ is large and the corresponding universal cover
  contains no flat planes.  Figure \ref{lotlotlink} represents the link $X_1$.  As we can see, there are no new circuits
of length less than five created in the link when the edges
corresponding to the extra square 2-cell  are added.  Thus, $X_1$
is non-positively curved and the poison corners from each original
complex are still present in $X_1$ prohibiting the development of
a flat plane in the universal cover, $\widetilde{X_1}$. Therefore
$\widetilde{X_1}$ is $\delta$-hyperbolic, and $G_1$ is a
hyperbolic group.  Theorem \ref{maps} implies our group is
free-by-cyclic and  there is a map from $G_1$ onto $\Z^2$, thus
completing the proof of Theorem \ref{main}.

   \end{proof}


   We leave it to the reader to check that $\langle L_{(k_1, a_i)}, L_{(k_2, b_i)}\hspace{0.2cm} | \hspace{0.4cm}a_0b_2a_1^{-1}b_0^{-1} \rangle$ is also free-by-cyclic and hyperbolic for all $k_1, k_2  \geq 4$.
    The following lemma gives us a formula for calculating the free rank of the kernel when there is an epimorphism from these groups onto $\Z$.


 \begin{figure}[htbp] \begin{center} 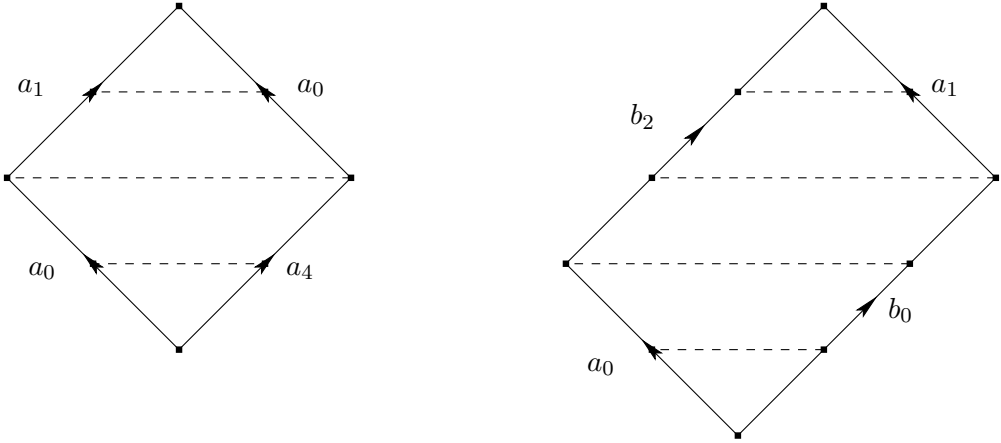 \caption{Example 2-cells of $K$ for the case $m=2$ and $n=3$, with the edges sub-divided.  The dashed lines represent the  pre-image of $p$ in these 2-cells.} \label{growthcell} \end{center} \end{figure}


   \begin{lemma}
   \label{growth}

Let  $K$ be a 2-complex constructed as in Theorem \ref{maps} such that
$K_j = X_{L_{(a_i,k_j)}}$ for $k_1, k_2 \geq 4$ and the added 2-cell reads $ a_0^{\pm1}b_2^{\pm1}a_1^{\mp1}b_0^{\mp1}$ around the boundary.
  Let $f_j:K_j \rightarrow S^1$, $j=1,2$, be circled-valued Morse functions such that  $m =|{f_1}_*(a_i)|$ and $n=|{f_2}_*(b_i)|$, and let $g:K\rightarrow S^1$ be the circle valued Morse function obtained from the construction in Theorem \ref{maps}.  Then, if $m, n \neq 0$ and $gcd(m, n) = 1$,  $\chi(g^{-1}(p)) =-(k_1m+k_2n)$.

     \end{lemma}

  \begin{proof}  First we will examine the point pre-images of the map $f_j:K_j
\rightarrow S^1$. Since we are dealing with affine cell-complexes
and affine maps, these pre-images are graphs.
 Since $K_j = X_{L_{(a_i, k)}}$,  $K_j$ is made up of $k$ square 2-cells with $(k+1)$ different edge
labels.
  For $|f_{j_*}(a_i)| = m$, the pre-image of the basepoint $p$  of $S^1$ is the union of its pre-images in the
2-cells.  Each
  edge will have to be subdivided which will add $(m-1)$ new vertices for each edge label (i.e. $(k+1)(m-1)
+1$
  total vertices).  The number of edges in the pre-image in a single 2-cell is equal to the number
of pairs of
  opposing vertices in the 2-cell, which is $2(m-1) + 1$.
  The left-hand side of Figure \ref{growthcell} shows a sample cell for the case $m=2$. The resulting graph has Euler characteristic:
$$(k+1)(m-1) +1 -
  k(2(m-1) +1)= -(k-1)m.$$

 Now we will look at the pre-image in the added 2-cell.  (The case where $m=2$ and $n=3$ is shown in Figure \ref{growthcell}.)  We will not count vertices since they
are counted in the building block complexes.  The number of edges
is $(m+n-1)$ which is again the number of opposing pairs of vertices in the cell.

Since the 2-complex $K$ is composed of $X_{L_{(k_1, a_i)}}$,
$X_{L_{(k_2, b_i)}}$ and the added 2-cell, after making the correction for the over count
of the original vertex by subtracting 1, we get:  $$\chi(g^{-1}(p)) = -(k_1-1)m
+ -(k_2 -1)n  - (m+n-1) - 1 = -(k_1m + k_2n).$$
\end{proof}

\begin{remark}
For $gcd(m,n)= d >1$ let $m'=m/d, n'=n/d$; then
$gcd(m',n')=1$.  Using lemma \ref{growth} we have
$\chi(h^{-1}(p)) = 1 -d(k_1m'+k_2n')-d$ where $h$ is the
composition of the circle-valued Morse function $g$ in the lemma
above with a  map which is multiplication by $d$ at the fundamental
group level i.e,
$h_* = (\times d)\circ g_*$.
\end{remark}

\begin{remark}
From the above result we see that the free rank of any group
constructed in this manner is $(k_1m + k_2n +1)$.  The minimal free rank
of $G_1$ above is therefore 9, and using these methods, we can
construct a group with any minimal free rank $k\geq 9$.
\end{remark}

Finally we claim that the examples we obtained in this section so
far have reducible automorphisms associated with them.
 Bestvina and Handel \cite{BeHa} gave the following definition of reducible and
 irreducible automorphisms.

\begin{definition}
     An automorphism of a group G is \emph{reducible} if it preserves the conjugacy class of a proper free
      factor of G.  If the automorphism is not reducible then it is \emph{irreducible}.

\end{definition}

The main example produced in this section is the group
$F_9\rtimes\Z$. This group was constructed as a one-relator
quotient of a free product of two free-by-cyclic groups (both of
free rank 4). We will show that there is a reducible automorphism
associated with this group.


 \begin{figure}[htbp] \begin{center} 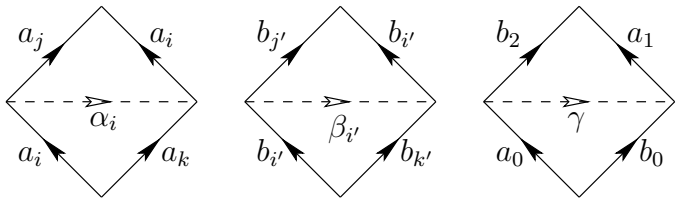 \caption{The generators of $F_9$.  The generators are labeled according to the edge labels of the corresponding 2-cell which appear twice in the cell.} \label{reducible} \end{center} \end{figure}


\begin{proposition}
\label{reducible} There exists a reducible automorphism associated
with the group $F_9\rtimes\Z$ constructed in the example above.
\end{proposition}

\begin{proof}

 Consider the pre-image of the
base-point of the target circle
 under the Morse function $g: X \rightarrow S^1$ (from Theorem \ref{maps}). This pre-image
 can be described as a union of the pre-images under $f_1:X_{L_{(a_i)}} \rightarrow S^1$ and
 $f_2:X_{L_{(b_i)}} \rightarrow S^1$ along with the pre-image in the extra 2-cell .
  The pre-images under $f_1$ and  $f_2$ are bouquets of 4 circles.
  Let $\alpha_0,.....,\alpha_3$ be the generators of the
 free group (denoted by $F_{\alpha_i}$) corresponding to the pre-image under
 $f_1$. Similarly let $\beta_0,....,\beta_3$ be the generators of the free group (denoted by
$F_{\beta_i}$) corresponding to the pre-image
 under $f_2$. Also, let the generator corresponding to the pre-image in the
 remaining 2-cell be $\gamma$.  Figure \ref{reducible} shows generic square
 2-cells of the 2-complex $X$ and the edge running through the middle of each square 2-cell representing the generators $\alpha_i, \beta_{i'}$ and $\gamma$.


\begin{figure}[htbp] \begin{center} 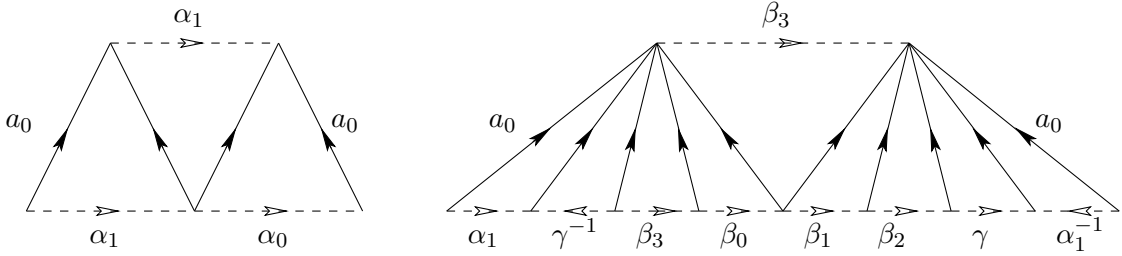 \caption{The effect of conjugating by $a_0$ on the generators $\alpha_1$ and $\beta_3$.  The automorphism is given by reading the labels on the bottom of the figures.} \label{autos} \end{center} \end{figure}


 If we conjugate by a generator from the
group $L_{(a_i)}$, say $a_0$, then $F_{\alpha_i}$ is invariant.
This can be shown by examining what effect this conjugation has on the generators.
$$ a_0 \alpha_0
{a_0}^{-1} = \alpha_3^{-1} \alpha_2^{-1} \alpha_1^{-1},  \hspace{0.3cm} a_0
\alpha_1 {a_0}^{-1} = \alpha_1 {\alpha_0},$$
$$a_0
\alpha_2 {a_0}^{-1} = \alpha_1 \alpha_2 {\alpha_0},
\hspace{0.3cm} a_0 \alpha_3 {a_0}^{-1} = \alpha_1 \alpha_2
\alpha_3 \alpha_0 $$
 Figure \ref{autos} illustrates how these relations are obtained.
 The automorphism
corresponding to the above conjugation actually keeps $F_{\alpha_i}$ invariant.  Thus this automorphism is reducible. We
can also see that $F_{\beta_i}$ is not invariant, but it is invariant up to conjugation.   Conjugating by $a_0$, gives us $\alpha \gamma^{-1} F_{\beta_i} (\alpha
\gamma^{-1})^{-1}$ as seen below.
$$a_0 \beta_0 a_0^{-1} = \alpha_1 \gamma^{-1} \beta_2^{-1} \beta_1^{-1}
\beta_3^{-1} \gamma \alpha_1^{-1}, \hspace{0.3cm} a_0 \beta_1 a_0^{-1}=
\alpha_1 \gamma^{-1} \beta_2^{-1} \beta_0 \beta_1 \beta_2 \gamma
\alpha_1^{-1}$$
$$a_0 \beta_2 a_0^{-1}= \alpha_1 \gamma^{-1} \beta_0 \beta_1
\beta_2 \gamma \alpha_1^{-1},\hspace{0.3cm} a_0 \beta_3 a_0^{-1} =
\alpha_1\gamma^{-1} \beta_3 \beta_0 \beta_1 \beta_2 \gamma \alpha_1^{-1}$$
In fact, conjugating by any $a_i$ generator will keep $F_{\alpha_i}$ invariant just as conjugating by any $b_i$ generator will keep $F_{\beta_i}$ invariant.  If we conjugate by other elements that map to $1\in \Z$, we will still preserve the conjugacy classes of these factors.  For example:
$$ (b_0 a_1^{-1}a_0) F_{\beta_i}( a_0^{-1} a_1 b_0^{-1}) \cong (b_0a_1^{-1})\alpha \gamma^{-1} F_{\beta_i}\gamma \alpha^{-1}(a_1b_0^{-1})$$
Since $b_0a_0^{-1}$ is in the kernel of $g_*$, we see that the conjugacy class of $F_{\beta_i}$ is still preserved under conjugation in $F_9 \rtimes \Z$.  Thus any automorphism will be reducible.

\end{proof}

This same technique can be used to give reducible automorphisms  for any group of the form, $\langle L_{(k_1, a_i)}, L_{(k_2,
b_i)}\hspace{0.2cm} | \hspace{0.4cm}a_0b_2a_1^{-1}b_0^{-1}
\rangle$, for $k_1, k_2 \geq 4$.


\section{A Labeled Oriented Forest Example }
By construction, the examples produced in the previous section are necessarily reducible.  In this section we will use a different method of construction to produce an example that does not immediately appear to be reducible.  This construction is similar to that presented in Section 4, but instead of taking a free product of two separate groups and adding a relation, we will simply add a relation to an already existing group.
The group we will begin with will be defined from a Labeled Oriented Forest (LOF).


\begin{figure}[htbp] \begin{center} 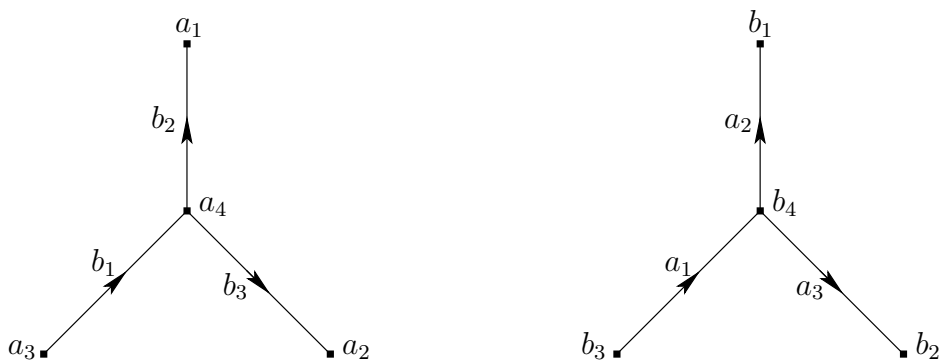 \caption{The LOF corresponding to the group $G_F$.} \label{lof} \end{center} \end{figure}


\label{LOF1}
\begin{example}
\label{LOF} (Labeled Oriented Forest):  Our LOF will consist of
the two trees shown in Figure \ref{lof}.  The first tree has
vertices labeled by $a_1, . . ., a_4$ and the second has vertices
labeled by $b_1, . . ., b_4$.  The edges of each tree are labeled by
vertices of the opposite tree.  The presentation 2-complex $X$ of
this LOF, is comprised of six square 2-cells and admits a
circle-valued Morse function $f: X\rightarrow S^1$.   The ascending
and descending links are 2-component forests.
 The corresponding group is: $$G_F = \langle a_i, b_i, i=1..4\hspace{0.3cm}| \hspace{0.3cm}
a_1^{b_2} = a_4,\hspace{0.1cm} a_2^{b_3} = a_4,
\hspace{0.1cm}a_4^{b_1} = a_3,\hspace{0.1cm} b_1^{a_2} =
b_4,\hspace{0.1cm} b_2^{a_3} =b_4,\hspace{0.1cm} b_4^{a_1} = b_3
\rangle . $$

The link of the complex (see Figure \ref{loflinks}) has no bigons or
empty triangles, therefore the link is large.  No flat planes
develop in the universal cover because every 2-cell of our complex
has not just one, but two poison corners.

We complete the
construction by adding an extra square 2-cell that represents the relation $a_4b_1 =
b_4a_1$.  The Morse function $f$ extends to a map $g_2$ over
the new complex which we
 will call $X_2$.  The corresponding group $G_2$ is denoted by:  $\langle G_F\hspace{0.3cm}|\hspace{0.3cm} a_4b_1a_1^{-1}b_4^{-1}\rangle$.

 The ascending and descending links corresponding to the new complex are now trees as the new
 relation contributes one edge each to both links joining up the 2-components of the forests.  Hence, the corresponding group is free-by-cyclic by Proposition \ref{FbyC}.
   The addition of this cell has added length 4 circuits to our link, so we must now check the new link to show
   that $G_2$ is hyperbolic.  The new link is shown on the right  in Figure
   \ref{loflinks}. The four edges contributed by the new
   2-cell are denoted here by dashed-dotted lines between the following vertices:
   ${b_4}^+{a_1}^-$, ${b_4}^-{a_4}^-$,  ${a_4}^+{b_1}^-$, and
   ${a_1}^+{b_1}^+$.
   We can see that where we previously had twelve poison corners, we now have only  four.
    Consequently,  our new 2-cell, and two of our original 2-cells no longer have poison corners.
    The following lemma will show that $G_2$ is indeed hyperbolic.


\begin{figure}[htbp] \begin{center} 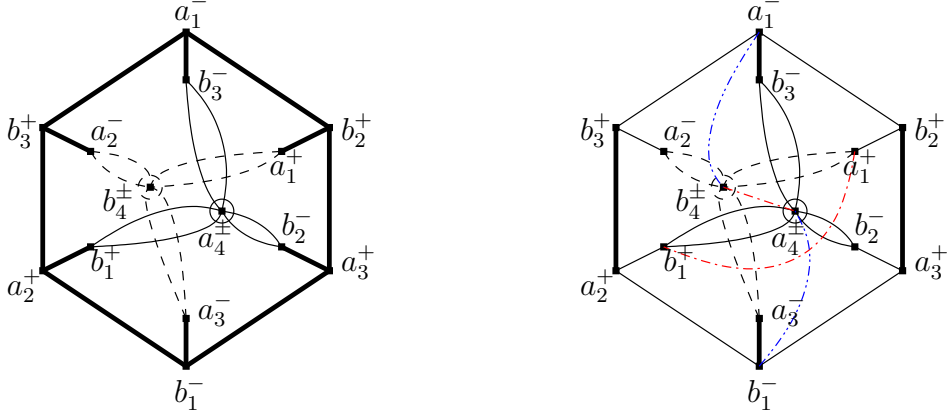 \caption{$Link(X, v)$ and $Link(X_2, v)$.  The circles indicate suspension points as in Figure \ref{lotandlink}.  Poison Corners are indicated by bold edges.  Note that although the one dot-dash edges and the two dot-dash edges appear to meet at the suspension points, they do not.  A careful check shows that there are no circuits of length less than four. } \label{loflinks} \end{center} \end{figure}


\begin{lemma}
   $G_2$ is hyperbolic.
\end{lemma}

 \begin{proof}  In order to prove this lemma, we must check that the universal cover of $X_2$ is $\delta$-hyperbolic,
 i.e., we must check that no flat planes can develop in the universal cover $\tilde{X_2}$ .
 Any such plane would have to be made up of the new square 2-cell, and only those original 2-cells which now have no poison
 corners.  We will attempt to build a plane starting with the new cell and show
 that it is impossible to develop a flat plane, thus proving the lemma.
 The 2-cells available to us are shown at the top of Figure \ref{poison}.
 In the bottom of Figure \ref{poison}, we see that it does not take long to discover a vertex around which we cannot find a Euclidean neighborhood
 (indicated by the large arrow).

 \end{proof}


\begin{figure}[htbp] \begin{center} 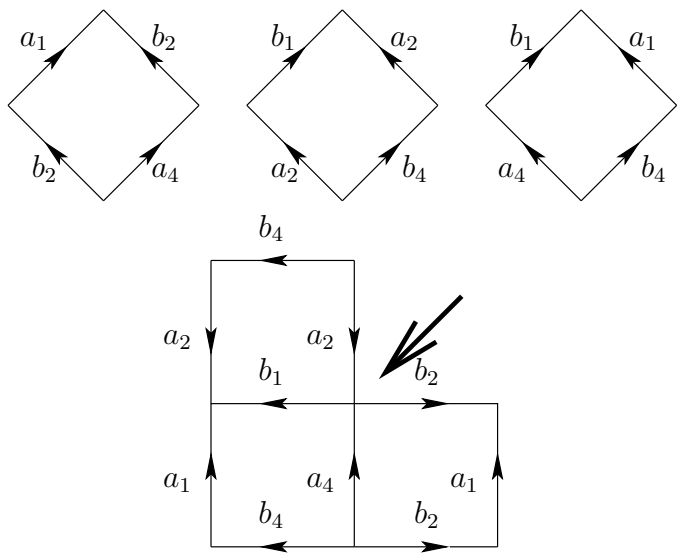 \caption{The 2-cells of $X_2$ which have no poison corners fail to develop into a plane in the universal cover.} \label{poison}\end{center} \end{figure}


 $G_2$ can be mapped onto $\Z^2$ because the $a_i$ generators and $b_i$ generators can be sent independently to generators of $\Z^2$.  Hence, $G_2$ fibers in infinitely many ways over $\Z$.
 The following lemma investigates the ranks of the kernels of these maps.
  \begin {lemma}
   Let $G_2$ and $g_2$ be as in Example \ref{LOF}.  Suppose that $m = {g_2}_*(a_i)$ and $n={g_2}_\ast(b_i)$ for $i=1,. .,4$.  Then if $m, n \neq 0$ and $gcd(m,n)= 1$, the rank of the kernel of the map ${g_2}_\ast : G_2 \rightarrow Z$ is given by the following expression:  $1+(3|m|+3|n|)$.
\end{lemma}

 \begin{proof}  By construction,  $m>0$ implies $n>0$, so we will assume without loss of generality that $m,n > 0$. In order
to find the rank,
  we will examine the Euler characteristic of the graph which is the pre-image of the base-point
of
  $S^1$ under the map $g_2$.   Since the 2-complex contains seven 2-cells, we will
  look at what the effect of changing the map $g_{2_\ast} $ has on each 2-cell.  In the presentation 2-complex, all the 2-cells look the same.  For $m=1$ and $n=1$, the pre-image has one vertex and seven edges (one from each 2-cell), so
the
   rank of the kernel is 7.
   Otherwise we must subdivide the edges, adding new vertices.   The number of new vertices corresponding to
each
   $a_i$, $b_i$ edge is  $4(m-1),4(n-1)$ respectively.   The number of edges in each 2-cell is $(m+n-1)$  (see the right-hand side Figure \ref{growthcell}),  so in all
there are $7(m+n-1)$ edges.  Counting the original vertex, the
total number of vertices is $4m+4n-7$, so the  Euler characteristic
of the graph is $-(3m+3n)$.  Therefore the rank of the free kernel is
$1+3m+3n$.

   \end{proof}

   \end{example}

   \begin{remark}For the example constructed in this section, we can show that there is a
  topological representative of an outer automorphism which is irreducible.  We will show this using definitions and techniques from Bestvina-Handel's paper \cite{BeHa}, Section 1.
  
  Let $R_n$ be the bouquet on $n$-circles.  A marked graph $G$ is a
graph with a homotopy equivalence $\tau: R_n \rightarrow G$.  A
homotopy equivalence $f: G \rightarrow G $ determines an outer
automorphism $ \mathcal{O}$, of $F_n$. If $\textit{V}$ is the
vertex set of $G$ and $f|_{ (G\setminus \textit{V})}$ is locally
injective, then we say that $f: G \rightarrow G$ is a
\emph{topological representative} of $ \mathcal{O}$.
 The transition matrix $M$ for a topological representative $f$ has entries $a_{ij}$ defined as the
number of times that the $f$-image of the $j ^{th}$ edge crosses
the $i^{th}$ edge in either direction.  A nonnegative integral matrix $M$ is irreducible if for all $1
\leq i, j \leq dim(M)$, there exists $N(i,j)>0$ so that the
$ij^{th}$ entry of $M^{N(i,j)}$ is positive.  A topological representative, $f$, is irreducible if and
only if its transition matrix is irreducible.  An outer automorphism $\mathcal{O}$ is irreducible if and only if every topological representative of $\mathcal{O}$ is irreducible.

Consider  the map ${g_2}_*: G_2 \rightarrow \Z$ such that $1 =
{g_2}_*(a_i), 1 = {g_2}_*(b_i)$.  Using the same
naming convention as in Proposition \ref{reducible}, let the generators of the
Ker(${g_2}_*$) be $\alpha_1, \alpha_2, \alpha_3, \beta_1, \beta_2,
\beta_3, \gamma$.  If we conjugate the generators by $a_3$ say, then we have the following
expressions:
$$ a_3 \alpha_1{a_3}^{-1} = \alpha_3^{-1} \beta_2^{-1} \beta_3^{-1}
\alpha_2^{-1}\gamma \beta_2 \alpha_3,
 \hspace{0.4cm} a_3
\alpha_2 {a_3}^{-1} = \alpha_3^{-1} \beta_2^{-1} \gamma^{-1} \alpha_2
\alpha_1\beta_2\alpha_3,$$
$$ a_3 \alpha_3 {a_3}^{-1} = \alpha_1 \beta_2 \alpha_3,
 \hspace{0.35cm}a_3 \beta_1
{a_3}^{-1} = \alpha_3^{-1} \beta_2^{-1} \gamma^{-1},
\hspace{0.35cm} a_3
\beta_2 {a_3}^{-1} = \alpha_3^{-1} \beta_1 \gamma \beta_2
\alpha_3,$$
$$a_3\beta_3{a_3}^{-1} = \alpha_3^{-1} \beta_2^{-1} \gamma^{-1}
\alpha_2 \beta_3 \beta_1 \gamma \beta_2
\alpha_3,
\hspace{0.4cm} a_3 \gamma {a_3}^{-1} = \alpha_3^{-1}
\beta_2^{-1} \gamma^{-1} \beta_1^{-1} \alpha_1\beta_2 \alpha_3$$
This gives us a topological representative of an outer
automorphism $\mathcal{O}$ and the corresponding transition matrix
$M$ is a $7 \times 7$ nonnegative matrix.  Using applicable software, it is easy to calculate the powers of the is matrix and see that
all entries of $M^3$ are positive.  Therefore, we have an
irreducible topological representative of the outer automorphism; however, using this method it is not possible to check that every topological
representative is irreducible.
Gersten and Stallings \cite{GeSt} provide a method for showing that an outer automorphism is irreducible, using super irreducible characteristic polynomials of $\mathcal{O}_{ab}$.  This method also failed to give a positive result due to the abelianization process involved.

   \end{remark}


   \section{Questions}
Following are some open questions concerning hyperbolic, free-by-cyclic groups that fiber in infinitely many ways.  
 
 \emph{Question 1.}  We know that the minimal rank of
free kernels for hyperbolic, free-by-cyclic groups is 3.  Using our methods we can generate groups that fiber in infinitely many ways with free kernel of minimal rank $7$ and
any minimal rank $k \geq 9$.  Can we
find lower rank examples that fiber in infinitely many ways?

The impetus for this research was work done in the 3-manifold setting.  The next question relates to results obtained in that setting.

\emph{Question 2.} Every hyperbolic 3-manifold that fibers in infinitely many ways has an irreducible outer automorphism.  Is the automorphism in the LOF example
 irreducible?  Are there any other examples with irreducible automorphisms?

All of our examples have maps onto $\Z^2$.  A group is said to be large if a finite index subgroup maps onto $F_2$.  Our next question deals with possible connections to large groups.

\emph{Question 3.}  Are there free-by-cyclic groups that fiber in only one way, such as $L_{(a_i)}$, which have finite index subgroups which fiber in infinitely
many ways?  Will these groups have virtual first Betti number $> 1$.   Would these groups be large?

\bibliography{hyperfiber17}
\bibliographystyle{siam}

\end{document}

%% file: lotandlink.eps_t
\begin{picture}(0,0)%
\epsfig{file=lotandlink.eps}%
\end{picture}%
\setlength{\unitlength}{3552sp}%
\begingroup\makeatletter\ifx\SetFigFont\undefined%
\gdef\SetFigFont#1#2#3#4#5{%
  \reset@font\fontsize{#1}{#2pt}%
  \fontfamily{#3}\fontseries{#4}\fontshape{#5}%
  \selectfont}%
\fi\endgroup%
\begin{picture}(7751,2905)(976,-3035)
\put(2926,-1411){\makebox(0,0)[lb]{\smash{{\SetFigFont{11}{13.2}{\familydefault}{\mddefault}{\updefault}{\color[rgb]{0,0,0}$a_2$}%
}}}}
\put(5926,-286){\makebox(0,0)[lb]{\smash{{\SetFigFont{11}{13.2}{\familydefault}{\mddefault}{\updefault}{\color[rgb]{0,0,0}$0^-$}%
}}}}
\put(8251,-286){\makebox(0,0)[lb]{\smash{{\SetFigFont{11}{13.2}{\familydefault}{\mddefault}{\updefault}{\color[rgb]{0,0,0}$1^+$}%
}}}}
\put(7801,-2011){\makebox(0,0)[lb]{\smash{{\SetFigFont{11}{13.2}{\familydefault}{\mddefault}{\updefault}{\color[rgb]{0,0,0}$1^-$}%
}}}}
\put(5926,-2986){\makebox(0,0)[lb]{\smash{{\SetFigFont{11}{13.2}{\familydefault}{\mddefault}{\updefault}{\color[rgb]{0,0,0}$3^+$}%
}}}}
\put(8251,-2986){\makebox(0,0)[lb]{\smash{{\SetFigFont{11}{13.2}{\familydefault}{\mddefault}{\updefault}{\color[rgb]{0,0,0}$2^+$}%
}}}}
\put(7126,-1336){\makebox(0,0)[lb]{\smash{{\SetFigFont{11}{13.2}{\familydefault}{\mddefault}{\updefault}{\color[rgb]{0,0,0}$4^{\pm}$}%
}}}}
\put(7801,-1186){\makebox(0,0)[lb]{\smash{{\SetFigFont{11}{13.2}{\familydefault}{\mddefault}{\updefault}{\color[rgb]{0,0,0}$0^+$}%
}}}}
\put(6451,-2011){\makebox(0,0)[lb]{\smash{{\SetFigFont{11}{13.2}{\familydefault}{\mddefault}{\updefault}{\color[rgb]{0,0,0}$2^-$}%
}}}}
\put(6451,-1186){\makebox(0,0)[lb]{\smash{{\SetFigFont{11}{13.2}{\familydefault}{\mddefault}{\updefault}{\color[rgb]{0,0,0}$3^-$}%
}}}}
\put(2476,-1711){\makebox(0,0)[lb]{\smash{{\SetFigFont{11}{13.2}{\familydefault}{\mddefault}{\updefault}{\color[rgb]{0,0,0}$a_4$}%
}}}}
\put(3676,-1636){\makebox(0,0)[lb]{\smash{{\SetFigFont{11}{13.2}{\familydefault}{\mddefault}{\updefault}{\color[rgb]{0,0,0}$a_1$}%
}}}}
\put(2326,-286){\makebox(0,0)[lb]{\smash{{\SetFigFont{11}{13.2}{\familydefault}{\mddefault}{\updefault}{\color[rgb]{0,0,0}$a_0$}%
}}}}
\put(2326,-2911){\makebox(0,0)[lb]{\smash{{\SetFigFont{11}{13.2}{\familydefault}{\mddefault}{\updefault}{\color[rgb]{0,0,0}$a_2$}%
}}}}
\put(976,-1636){\makebox(0,0)[lb]{\smash{{\SetFigFont{11}{13.2}{\familydefault}{\mddefault}{\updefault}{\color[rgb]{0,0,0}$a_3$}%
}}}}
\put(2476,-1036){\makebox(0,0)[lb]{\smash{{\SetFigFont{11}{13.2}{\familydefault}{\mddefault}{\updefault}{\color[rgb]{0,0,0}$a_1$}%
}}}}
\put(1726,-1711){\makebox(0,0)[lb]{\smash{{\SetFigFont{11}{13.2}{\familydefault}{\mddefault}{\updefault}{\color[rgb]{0,0,0}$a_0$}%
}}}}
\put(2176,-2086){\makebox(0,0)[lb]{\smash{{\SetFigFont{11}{13.2}{\familydefault}{\mddefault}{\updefault}{\color[rgb]{0,0,0}$a_3$}%
}}}}
\end{picture}%

%% file: lotlotlink.eps_t
\begin{picture}(0,0)%
\epsfig{file=lotlotlink.eps}%
\end{picture}%
\setlength{\unitlength}{3552sp}%
\begingroup\makeatletter\ifx\SetFigFont\undefined%
\gdef\SetFigFont#1#2#3#4#5{%
  \reset@font\fontsize{#1}{#2pt}%
  \fontfamily{#3}\fontseries{#4}\fontshape{#5}%
  \selectfont}%
\fi\endgroup%
\begin{picture}(7601,3441)(1126,-3208)
\put(7126, 89){\makebox(0,0)[lb]{\smash{{\SetFigFont{12}{14.4}{\familydefault}{\mddefault}{\updefault}{\color[rgb]{0,0,0}$b$}%
}}}}
\put(1126,-286){\makebox(0,0)[lb]{\smash{{\SetFigFont{11}{13.2}{\familydefault}{\mddefault}{\updefault}{\color[rgb]{0,0,0}$0^-$}%
}}}}
\put(3451,-286){\makebox(0,0)[lb]{\smash{{\SetFigFont{11}{13.2}{\familydefault}{\mddefault}{\updefault}{\color[rgb]{0,0,0}$1^+$}%
}}}}
\put(3001,-2011){\makebox(0,0)[lb]{\smash{{\SetFigFont{11}{13.2}{\familydefault}{\mddefault}{\updefault}{\color[rgb]{0,0,0}$1^-$}%
}}}}
\put(1126,-2986){\makebox(0,0)[lb]{\smash{{\SetFigFont{11}{13.2}{\familydefault}{\mddefault}{\updefault}{\color[rgb]{0,0,0}$3^+$}%
}}}}
\put(3451,-2986){\makebox(0,0)[lb]{\smash{{\SetFigFont{11}{13.2}{\familydefault}{\mddefault}{\updefault}{\color[rgb]{0,0,0}$2^+$}%
}}}}
\put(2326,-1336){\makebox(0,0)[lb]{\smash{{\SetFigFont{11}{13.2}{\familydefault}{\mddefault}{\updefault}{\color[rgb]{0,0,0}$4^{\pm}$}%
}}}}
\put(3001,-1186){\makebox(0,0)[lb]{\smash{{\SetFigFont{11}{13.2}{\familydefault}{\mddefault}{\updefault}{\color[rgb]{0,0,0}$0^+$}%
}}}}
\put(1651,-2011){\makebox(0,0)[lb]{\smash{{\SetFigFont{11}{13.2}{\familydefault}{\mddefault}{\updefault}{\color[rgb]{0,0,0}$2^-$}%
}}}}
\put(1651,-1186){\makebox(0,0)[lb]{\smash{{\SetFigFont{11}{13.2}{\familydefault}{\mddefault}{\updefault}{\color[rgb]{0,0,0}$3^-$}%
}}}}
\put(5926,-286){\makebox(0,0)[lb]{\smash{{\SetFigFont{11}{13.2}{\familydefault}{\mddefault}{\updefault}{\color[rgb]{0,0,0}$1^+$}%
}}}}
\put(8251,-286){\makebox(0,0)[lb]{\smash{{\SetFigFont{11}{13.2}{\familydefault}{\mddefault}{\updefault}{\color[rgb]{0,0,0}$0^-$}%
}}}}
\put(8251,-2986){\makebox(0,0)[lb]{\smash{{\SetFigFont{11}{13.2}{\familydefault}{\mddefault}{\updefault}{\color[rgb]{0,0,0}$3^+$}%
}}}}
\put(5926,-2986){\makebox(0,0)[lb]{\smash{{\SetFigFont{11}{13.2}{\familydefault}{\mddefault}{\updefault}{\color[rgb]{0,0,0}$2^+$}%
}}}}
\put(6451,-1186){\makebox(0,0)[lb]{\smash{{\SetFigFont{11}{13.2}{\familydefault}{\mddefault}{\updefault}{\color[rgb]{0,0,0}$0^+$}%
}}}}
\put(7126,-1336){\makebox(0,0)[lb]{\smash{{\SetFigFont{11}{13.2}{\familydefault}{\mddefault}{\updefault}{\color[rgb]{0,0,0}$4^{\pm}$}%
}}}}
\put(7801,-1186){\makebox(0,0)[lb]{\smash{{\SetFigFont{11}{13.2}{\familydefault}{\mddefault}{\updefault}{\color[rgb]{0,0,0}$3^-$}%
}}}}
\put(7801,-2011){\makebox(0,0)[lb]{\smash{{\SetFigFont{11}{13.2}{\familydefault}{\mddefault}{\updefault}{\color[rgb]{0,0,0}$2^-$}%
}}}}
\put(6451,-2011){\makebox(0,0)[lb]{\smash{{\SetFigFont{11}{13.2}{\familydefault}{\mddefault}{\updefault}{\color[rgb]{0,0,0}$1^-$}%
}}}}
\put(2326, 89){\makebox(0,0)[lb]{\smash{{\SetFigFont{12}{14.4}{\familydefault}{\mddefault}{\updefault}{\color[rgb]{0,0,0}$a$}%
}}}}
\end{picture}%

%% file: growthcell.eps_t
\begin{picture}(0,0)%
\epsfig{file=growthcell.eps}%
\end{picture}%
\setlength{\unitlength}{3552sp}%
\begingroup\makeatletter\ifx\SetFigFont\undefined%
\gdef\SetFigFont#1#2#3#4#5{%
  \reset@font\fontsize{#1}{#2pt}%
  \fontfamily{#3}\fontseries{#4}\fontshape{#5}%
  \selectfont}%
\fi\endgroup%
\begin{picture}(6988,3088)(1757,-2805)
\put(3751,-1636){\makebox(0,0)[lb]{\smash{{\SetFigFont{11}{13.2}{\familydefault}{\mddefault}{\updefault}{\color[rgb]{0,0,0}$a_4$}%
}}}}
\put(5851,-2311){\makebox(0,0)[lb]{\smash{{\SetFigFont{11}{13.2}{\familydefault}{\mddefault}{\updefault}{\color[rgb]{0,0,0}$a_0$}%
}}}}
\put(6151,-586){\makebox(0,0)[lb]{\smash{{\SetFigFont{11}{13.2}{\familydefault}{\mddefault}{\updefault}{\color[rgb]{0,0,0}$b_2$}%
}}}}
\put(8251,-361){\makebox(0,0)[lb]{\smash{{\SetFigFont{11}{13.2}{\familydefault}{\mddefault}{\updefault}{\color[rgb]{0,0,0}$a_1$}%
}}}}
\put(7951,-1936){\makebox(0,0)[lb]{\smash{{\SetFigFont{11}{13.2}{\familydefault}{\mddefault}{\updefault}{\color[rgb]{0,0,0}$b_0$}%
}}}}
\put(1951,-1636){\makebox(0,0)[lb]{\smash{{\SetFigFont{11}{13.2}{\familydefault}{\mddefault}{\updefault}{\color[rgb]{0,0,0}$a_0$}%
}}}}
\put(1876,-361){\makebox(0,0)[lb]{\smash{{\SetFigFont{11}{13.2}{\familydefault}{\mddefault}{\updefault}{\color[rgb]{0,0,0}$a_1$}%
}}}}
\put(3826,-361){\makebox(0,0)[lb]{\smash{{\SetFigFont{11}{13.2}{\familydefault}{\mddefault}{\updefault}{\color[rgb]{0,0,0}$a_0$}%
}}}}
\end{picture}%

%% file: reducible.eps_t
\begin{picture}(0,0)%
\epsfig{file=reducible.eps}%
\end{picture}%
\setlength{\unitlength}{3947sp}%
\begingroup\makeatletter\ifx\SetFigFont\undefined%
\gdef\SetFigFont#1#2#3#4#5{%
  \reset@font\fontsize{#1}{#2pt}%
  \fontfamily{#3}\fontseries{#4}\fontshape{#5}%
  \selectfont}%
\fi\endgroup%
\begin{picture}(4448,1224)(1189,-1573)
\put(3226,-1186){\makebox(0,0)[lb]{\smash{{\SetFigFont{12}{14.4}{\familydefault}{\mddefault}{\updefault}{\color[rgb]{0,0,0}$\beta_{i'}$}%
}}}}
\put(1276,-586){\makebox(0,0)[lb]{\smash{{\SetFigFont{12}{14.4}{\familydefault}{\mddefault}{\updefault}{\color[rgb]{0,0,0}$a_j$}%
}}}}
\put(2101,-586){\makebox(0,0)[lb]{\smash{{\SetFigFont{12}{14.4}{\familydefault}{\mddefault}{\updefault}{\color[rgb]{0,0,0}$a_i$}%
}}}}
\put(2776,-586){\makebox(0,0)[lb]{\smash{{\SetFigFont{12}{14.4}{\familydefault}{\mddefault}{\updefault}{\color[rgb]{0,0,0}$b_{j'}$}%
}}}}
\put(2776,-1336){\makebox(0,0)[lb]{\smash{{\SetFigFont{12}{14.4}{\familydefault}{\mddefault}{\updefault}{\color[rgb]{0,0,0}$b_{i'}$}%
}}}}
\put(3601,-586){\makebox(0,0)[lb]{\smash{{\SetFigFont{12}{14.4}{\familydefault}{\mddefault}{\updefault}{\color[rgb]{0,0,0}$b_{i'}$}%
}}}}
\put(4276,-586){\makebox(0,0)[lb]{\smash{{\SetFigFont{12}{14.4}{\familydefault}{\mddefault}{\updefault}{\color[rgb]{0,0,0}$b_2$}%
}}}}
\put(5101,-586){\makebox(0,0)[lb]{\smash{{\SetFigFont{12}{14.4}{\familydefault}{\mddefault}{\updefault}{\color[rgb]{0,0,0}$a_1$}%
}}}}
\put(4276,-1336){\makebox(0,0)[lb]{\smash{{\SetFigFont{12}{14.4}{\familydefault}{\mddefault}{\updefault}{\color[rgb]{0,0,0}$a_0$}%
}}}}
\put(1276,-1336){\makebox(0,0)[lb]{\smash{{\SetFigFont{12}{14.4}{\familydefault}{\mddefault}{\updefault}{\color[rgb]{0,0,0}$a_i$}%
}}}}
\put(2176,-1336){\makebox(0,0)[lb]{\smash{{\SetFigFont{12}{14.4}{\familydefault}{\mddefault}{\updefault}{\color[rgb]{0,0,0}$a_k$}%
}}}}
\put(3676,-1336){\makebox(0,0)[lb]{\smash{{\SetFigFont{12}{14.4}{\familydefault}{\mddefault}{\updefault}{\color[rgb]{0,0,0}$b_{k'}$}%
}}}}
\put(5176,-1336){\makebox(0,0)[lb]{\smash{{\SetFigFont{12}{14.4}{\familydefault}{\mddefault}{\updefault}{\color[rgb]{0,0,0}$b_0$}%
}}}}
\put(1726,-1111){\makebox(0,0)[lb]{\smash{{\SetFigFont{12}{14.4}{\familydefault}{\mddefault}{\updefault}{\color[rgb]{0,0,0}$\alpha_i$}%
}}}}
\put(4726,-1111){\makebox(0,0)[lb]{\smash{{\SetFigFont{12}{14.4}{\familydefault}{\mddefault}{\updefault}{\color[rgb]{0,0,0}$\gamma$}%
}}}}
\end{picture}%

%% file: autos.eps_t
\begin{picture}(0,0)%
\epsfig{file=autos.eps}%
\end{picture}%
\setlength{\unitlength}{3473sp}%
\begingroup\makeatletter\ifx\SetFigFont\undefined%
\gdef\SetFigFont#1#2#3#4#5{%
  \reset@font\fontsize{#1}{#2pt}%
  \fontfamily{#3}\fontseries{#4}\fontshape{#5}%
  \selectfont}%
\fi\endgroup%
\begin{picture}(7962,1789)(2251,-3044)
\put(2851,-2986){\makebox(0,0)[lb]{\smash{{\SetFigFont{11}{13.2}{\familydefault}{\mddefault}{\updefault}{\color[rgb]{0,0,0}$\alpha_1$}%
}}}}
\put(4051,-2986){\makebox(0,0)[lb]{\smash{{\SetFigFont{11}{13.2}{\familydefault}{\mddefault}{\updefault}{\color[rgb]{0,0,0}$\alpha_0$}%
}}}}
\put(5701,-2161){\makebox(0,0)[lb]{\smash{{\SetFigFont{11}{13.2}{\familydefault}{\mddefault}{\updefault}{\color[rgb]{0,0,0}$a_0$}%
}}}}
\put(7951,-2986){\makebox(0,0)[lb]{\smash{{\SetFigFont{11}{13.2}{\familydefault}{\mddefault}{\updefault}{\color[rgb]{0,0,0}$\beta_1$}%
}}}}
\put(3451,-1411){\makebox(0,0)[lb]{\smash{{\SetFigFont{11}{13.2}{\familydefault}{\mddefault}{\updefault}{\color[rgb]{0,0,0}$\alpha_1$}%
}}}}
\put(2251,-2161){\makebox(0,0)[lb]{\smash{{\SetFigFont{11}{13.2}{\familydefault}{\mddefault}{\updefault}{\color[rgb]{0,0,0}$a_0$}%
}}}}
\put(4576,-2161){\makebox(0,0)[lb]{\smash{{\SetFigFont{11}{13.2}{\familydefault}{\mddefault}{\updefault}{\color[rgb]{0,0,0}$a_0$}%
}}}}
\put(9601,-2161){\makebox(0,0)[lb]{\smash{{\SetFigFont{11}{13.2}{\familydefault}{\mddefault}{\updefault}{\color[rgb]{0,0,0}$a_0$}%
}}}}
\put(5551,-2986){\makebox(0,0)[lb]{\smash{{\SetFigFont{11}{13.2}{\familydefault}{\mddefault}{\updefault}{\color[rgb]{0,0,0}$\alpha_1$}%
}}}}
\put(6151,-2986){\makebox(0,0)[lb]{\smash{{\SetFigFont{11}{13.2}{\familydefault}{\mddefault}{\updefault}{\color[rgb]{0,0,0}$\gamma^{-1}$}%
}}}}
\put(6751,-2986){\makebox(0,0)[lb]{\smash{{\SetFigFont{11}{13.2}{\familydefault}{\mddefault}{\updefault}{\color[rgb]{0,0,0}$\beta_3$}%
}}}}
\put(7351,-2986){\makebox(0,0)[lb]{\smash{{\SetFigFont{11}{13.2}{\familydefault}{\mddefault}{\updefault}{\color[rgb]{0,0,0}$\beta_0$}%
}}}}
\put(8476,-2986){\makebox(0,0)[lb]{\smash{{\SetFigFont{11}{13.2}{\familydefault}{\mddefault}{\updefault}{\color[rgb]{0,0,0}$\beta_2$}%
}}}}
\put(9151,-2986){\makebox(0,0)[lb]{\smash{{\SetFigFont{11}{13.2}{\familydefault}{\mddefault}{\updefault}{\color[rgb]{0,0,0}$\gamma$}%
}}}}
\put(9751,-2986){\makebox(0,0)[lb]{\smash{{\SetFigFont{11}{13.2}{\familydefault}{\mddefault}{\updefault}{\color[rgb]{0,0,0}$\alpha_1^{-1}$}%
}}}}
\put(7651,-1411){\makebox(0,0)[lb]{\smash{{\SetFigFont{11}{13.2}{\familydefault}{\mddefault}{\updefault}{\color[rgb]{0,0,0}$\beta_3$}%
}}}}
\end{picture}%

%% file: lof.eps_t
\begin{picture}(0,0)%
\epsfig{file=lof.eps}%
\end{picture}%
\setlength{\unitlength}{3947sp}%
\begingroup\makeatletter\ifx\SetFigFont\undefined%
\gdef\SetFigFont#1#2#3#4#5{%
  \reset@font\fontsize{#1}{#2pt}%
  \fontfamily{#3}\fontseries{#4}\fontshape{#5}%
  \selectfont}%
\fi\endgroup%
\begin{picture}(6161,2239)(1276,-2519)
\put(6976,-2461){\makebox(0,0)[lb]{\smash{{\SetFigFont{12}{14.4}{\familydefault}{\mddefault}{\updefault}{\color[rgb]{0,0,0}$b_2$}%
}}}}
\put(2326,-436){\makebox(0,0)[lb]{\smash{{\SetFigFont{12}{14.4}{\familydefault}{\mddefault}{\updefault}{\color[rgb]{0,0,0}$a_1$}%
}}}}
\put(2476,-1561){\makebox(0,0)[lb]{\smash{{\SetFigFont{12}{14.4}{\familydefault}{\mddefault}{\updefault}{\color[rgb]{0,0,0}$a_4$}%
}}}}
\put(6076,-1561){\makebox(0,0)[lb]{\smash{{\SetFigFont{12}{14.4}{\familydefault}{\mddefault}{\updefault}{\color[rgb]{0,0,0}$b_4$}%
}}}}
\put(1276,-2461){\makebox(0,0)[lb]{\smash{{\SetFigFont{12}{14.4}{\familydefault}{\mddefault}{\updefault}{\color[rgb]{0,0,0}$a_3$}%
}}}}
\put(2626,-2086){\makebox(0,0)[lb]{\smash{{\SetFigFont{12}{14.4}{\familydefault}{\mddefault}{\updefault}{\color[rgb]{0,0,0}$b_3$}%
}}}}
\put(2176,-1036){\makebox(0,0)[lb]{\smash{{\SetFigFont{12}{14.4}{\familydefault}{\mddefault}{\updefault}{\color[rgb]{0,0,0}$b_2$}%
}}}}
\put(1801,-1936){\makebox(0,0)[lb]{\smash{{\SetFigFont{12}{14.4}{\familydefault}{\mddefault}{\updefault}{\color[rgb]{0,0,0}$b_1$}%
}}}}
\put(3376,-2461){\makebox(0,0)[lb]{\smash{{\SetFigFont{12}{14.4}{\familydefault}{\mddefault}{\updefault}{\color[rgb]{0,0,0}$a_2$}%
}}}}
\put(5776,-1036){\makebox(0,0)[lb]{\smash{{\SetFigFont{12}{14.4}{\familydefault}{\mddefault}{\updefault}{\color[rgb]{0,0,0}$a_2$}%
}}}}
\put(5926,-436){\makebox(0,0)[lb]{\smash{{\SetFigFont{12}{14.4}{\familydefault}{\mddefault}{\updefault}{\color[rgb]{0,0,0}$b_1$}%
}}}}
\put(4876,-2461){\makebox(0,0)[lb]{\smash{{\SetFigFont{12}{14.4}{\familydefault}{\mddefault}{\updefault}{\color[rgb]{0,0,0}$b_3$}%
}}}}
\put(5401,-1936){\makebox(0,0)[lb]{\smash{{\SetFigFont{12}{14.4}{\familydefault}{\mddefault}{\updefault}{\color[rgb]{0,0,0}$a_1$}%
}}}}
\put(6226,-2086){\makebox(0,0)[lb]{\smash{{\SetFigFont{12}{14.4}{\familydefault}{\mddefault}{\updefault}{\color[rgb]{0,0,0}$a_3$}%
}}}}
\end{picture}%

%% file: loflinks.eps_t
\begin{picture}(0,0)%
\epsfig{file=loflinks.eps}%
\end{picture}%
\setlength{\unitlength}{3947sp}%
\begingroup\makeatletter\ifx\SetFigFont\undefined%
\gdef\SetFigFont#1#2#3#4#5{%
  \reset@font\fontsize{#1}{#2pt}%
  \fontfamily{#3}\fontseries{#4}\fontshape{#5}%
  \selectfont}%
\fi\endgroup%
\begin{picture}(6361,2614)(4876,-2744)
\put(9001,-1036){\makebox(0,0)[lb]{\smash{{\SetFigFont{12}{14.4}{\familydefault}{\mddefault}{\updefault}{\color[rgb]{0,0,0}$a_2^-$}%
}}}}
\put(5926,-286){\makebox(0,0)[lb]{\smash{{\SetFigFont{12}{14.4}{\familydefault}{\mddefault}{\updefault}{\color[rgb]{0,0,0}$a_1^-$}%
}}}}
\put(6976,-1861){\makebox(0,0)[lb]{\smash{{\SetFigFont{12}{14.4}{\familydefault}{\mddefault}{\updefault}{\color[rgb]{0,0,0}$a_3^+$}%
}}}}
\put(4876,-2011){\makebox(0,0)[lb]{\smash{{\SetFigFont{12}{14.4}{\familydefault}{\mddefault}{\updefault}{\color[rgb]{0,0,0}$a_2^+$}%
}}}}
\put(4876,-1036){\makebox(0,0)[lb]{\smash{{\SetFigFont{12}{14.4}{\familydefault}{\mddefault}{\updefault}{\color[rgb]{0,0,0}$b_3^+$}%
}}}}
\put(5926,-2686){\makebox(0,0)[lb]{\smash{{\SetFigFont{12}{14.4}{\familydefault}{\mddefault}{\updefault}{\color[rgb]{0,0,0}$b_1^-$}%
}}}}
\put(6976,-1036){\makebox(0,0)[lb]{\smash{{\SetFigFont{12}{14.4}{\familydefault}{\mddefault}{\updefault}{\color[rgb]{0,0,0}$b_2^+$}%
}}}}
\put(6076,-736){\makebox(0,0)[lb]{\smash{{\SetFigFont{12}{14.4}{\familydefault}{\mddefault}{\updefault}{\color[rgb]{0,0,0}$b_3^-$}%
}}}}
\put(6526,-1261){\makebox(0,0)[lb]{\smash{{\SetFigFont{12}{14.4}{\familydefault}{\mddefault}{\updefault}{\color[rgb]{0,0,0}$a_1^+$}%
}}}}
\put(6601,-1636){\makebox(0,0)[lb]{\smash{{\SetFigFont{12}{14.4}{\familydefault}{\mddefault}{\updefault}{\color[rgb]{0,0,0}$b_2^-$}%
}}}}
\put(6076,-1711){\makebox(0,0)[lb]{\smash{{\SetFigFont{12}{14.4}{\familydefault}{\mddefault}{\updefault}{\color[rgb]{0,0,0}$a_4^\pm$}%
}}}}
\put(5401,-1861){\makebox(0,0)[lb]{\smash{{\SetFigFont{12}{14.4}{\familydefault}{\mddefault}{\updefault}{\color[rgb]{0,0,0}$b_1^+$}%
}}}}
\put(5476,-1486){\makebox(0,0)[lb]{\smash{{\SetFigFont{12}{14.4}{\familydefault}{\mddefault}{\updefault}{\color[rgb]{0,0,0}$b_4^{\pm}$}%
}}}}
\put(6076,-2161){\makebox(0,0)[lb]{\smash{{\SetFigFont{12}{14.4}{\familydefault}{\mddefault}{\updefault}{\color[rgb]{0,0,0}$a_3^-$}%
}}}}
\put(5401,-1036){\makebox(0,0)[lb]{\smash{{\SetFigFont{12}{14.4}{\familydefault}{\mddefault}{\updefault}{\color[rgb]{0,0,0}$a_2^-$}%
}}}}
\put(9526,-286){\makebox(0,0)[lb]{\smash{{\SetFigFont{12}{14.4}{\familydefault}{\mddefault}{\updefault}{\color[rgb]{0,0,0}$a_1^-$}%
}}}}
\put(8476,-1036){\makebox(0,0)[lb]{\smash{{\SetFigFont{12}{14.4}{\familydefault}{\mddefault}{\updefault}{\color[rgb]{0,0,0}$b_3^+$}%
}}}}
\put(9526,-2686){\makebox(0,0)[lb]{\smash{{\SetFigFont{12}{14.4}{\familydefault}{\mddefault}{\updefault}{\color[rgb]{0,0,0}$b_1^-$}%
}}}}
\put(10576,-1036){\makebox(0,0)[lb]{\smash{{\SetFigFont{12}{14.4}{\familydefault}{\mddefault}{\updefault}{\color[rgb]{0,0,0}$b_2^+$}%
}}}}
\put(9676,-736){\makebox(0,0)[lb]{\smash{{\SetFigFont{12}{14.4}{\familydefault}{\mddefault}{\updefault}{\color[rgb]{0,0,0}$b_3^-$}%
}}}}
\put(10126,-1261){\makebox(0,0)[lb]{\smash{{\SetFigFont{12}{14.4}{\familydefault}{\mddefault}{\updefault}{\color[rgb]{0,0,0}$a_1^+$}%
}}}}
\put(10201,-1636){\makebox(0,0)[lb]{\smash{{\SetFigFont{12}{14.4}{\familydefault}{\mddefault}{\updefault}{\color[rgb]{0,0,0}$b_2^-$}%
}}}}
\put(9676,-1711){\makebox(0,0)[lb]{\smash{{\SetFigFont{12}{14.4}{\familydefault}{\mddefault}{\updefault}{\color[rgb]{0,0,0}$a_4^\pm$}%
}}}}
\put(9001,-1861){\makebox(0,0)[lb]{\smash{{\SetFigFont{12}{14.4}{\familydefault}{\mddefault}{\updefault}{\color[rgb]{0,0,0}$b_1^+$}%
}}}}
\put(9076,-1486){\makebox(0,0)[lb]{\smash{{\SetFigFont{12}{14.4}{\familydefault}{\mddefault}{\updefault}{\color[rgb]{0,0,0}$b_4^{\pm}$}%
}}}}
\put(9676,-2161){\makebox(0,0)[lb]{\smash{{\SetFigFont{12}{14.4}{\familydefault}{\mddefault}{\updefault}{\color[rgb]{0,0,0}$a_3^-$}%
}}}}
\put(10576,-1861){\makebox(0,0)[lb]{\smash{{\SetFigFont{12}{14.4}{\familydefault}{\mddefault}{\updefault}{\color[rgb]{0,0,0}$a_3^+$}%
}}}}
\put(8476,-2011){\makebox(0,0)[lb]{\smash{{\SetFigFont{12}{14.4}{\familydefault}{\mddefault}{\updefault}{\color[rgb]{0,0,0}$a_2^+$}%
}}}}
\end{picture}%

%% file: poison.eps_t
\begin{picture}(0,0)%
\epsfig{file=poison.eps}%
\end{picture}%
\setlength{\unitlength}{3947sp}%
\begingroup\makeatletter\ifx\SetFigFont\undefined%
\gdef\SetFigFont#1#2#3#4#5{%
  \reset@font\fontsize{#1}{#2pt}%
  \fontfamily{#3}\fontseries{#4}\fontshape{#5}%
  \selectfont}%
\fi\endgroup%
\begin{picture}(4448,3429)(1189,-3778)
\put(5176,-1411){\makebox(0,0)[lb]{\smash{{\SetFigFont{12}{14.4}{\familydefault}{\mddefault}{\updefault}{\color[rgb]{0,0,0}$b_4$}%
}}}}
\put(2776,-1786){\makebox(0,0)[lb]{\smash{{\SetFigFont{12}{14.4}{\familydefault}{\mddefault}{\updefault}{\color[rgb]{0,0,0}$b_4$}%
}}}}
\put(2176,-2461){\makebox(0,0)[lb]{\smash{{\SetFigFont{12}{14.4}{\familydefault}{\mddefault}{\updefault}{\color[rgb]{0,0,0}$a_2$}%
}}}}
\put(3076,-2461){\makebox(0,0)[lb]{\smash{{\SetFigFont{12}{14.4}{\familydefault}{\mddefault}{\updefault}{\color[rgb]{0,0,0}$a_2$}%
}}}}
\put(2776,-2686){\makebox(0,0)[lb]{\smash{{\SetFigFont{12}{14.4}{\familydefault}{\mddefault}{\updefault}{\color[rgb]{0,0,0}$b_1$}%
}}}}
\put(3751,-2686){\makebox(0,0)[lb]{\smash{{\SetFigFont{12}{14.4}{\familydefault}{\mddefault}{\updefault}{\color[rgb]{0,0,0}$b_2$}%
}}}}
\put(2176,-3361){\makebox(0,0)[lb]{\smash{{\SetFigFont{12}{14.4}{\familydefault}{\mddefault}{\updefault}{\color[rgb]{0,0,0}$a_1$}%
}}}}
\put(3076,-3361){\makebox(0,0)[lb]{\smash{{\SetFigFont{12}{14.4}{\familydefault}{\mddefault}{\updefault}{\color[rgb]{0,0,0}$a_4$}%
}}}}
\put(3976,-3361){\makebox(0,0)[lb]{\smash{{\SetFigFont{12}{14.4}{\familydefault}{\mddefault}{\updefault}{\color[rgb]{0,0,0}$a_1$}%
}}}}
\put(2776,-3586){\makebox(0,0)[lb]{\smash{{\SetFigFont{12}{14.4}{\familydefault}{\mddefault}{\updefault}{\color[rgb]{0,0,0}$b_4$}%
}}}}
\put(3751,-3586){\makebox(0,0)[lb]{\smash{{\SetFigFont{12}{14.4}{\familydefault}{\mddefault}{\updefault}{\color[rgb]{0,0,0}$b_2$}%
}}}}
\put(1276,-586){\makebox(0,0)[lb]{\smash{{\SetFigFont{12}{14.4}{\familydefault}{\mddefault}{\updefault}{\color[rgb]{0,0,0}$a_1$}%
}}}}
\put(2101,-586){\makebox(0,0)[lb]{\smash{{\SetFigFont{12}{14.4}{\familydefault}{\mddefault}{\updefault}{\color[rgb]{0,0,0}$b_2$}%
}}}}
\put(2101,-1411){\makebox(0,0)[lb]{\smash{{\SetFigFont{12}{14.4}{\familydefault}{\mddefault}{\updefault}{\color[rgb]{0,0,0}$a_4$}%
}}}}
\put(2851,-586){\makebox(0,0)[lb]{\smash{{\SetFigFont{12}{14.4}{\familydefault}{\mddefault}{\updefault}{\color[rgb]{0,0,0}$b_1$}%
}}}}
\put(4351,-586){\makebox(0,0)[lb]{\smash{{\SetFigFont{12}{14.4}{\familydefault}{\mddefault}{\updefault}{\color[rgb]{0,0,0}$b_1$}%
}}}}
\put(5101,-586){\makebox(0,0)[lb]{\smash{{\SetFigFont{12}{14.4}{\familydefault}{\mddefault}{\updefault}{\color[rgb]{0,0,0}$a_1$}%
}}}}
\put(2851,-1411){\makebox(0,0)[lb]{\smash{{\SetFigFont{12}{14.4}{\familydefault}{\mddefault}{\updefault}{\color[rgb]{0,0,0}$a_2$}%
}}}}
\put(3601,-586){\makebox(0,0)[lb]{\smash{{\SetFigFont{12}{14.4}{\familydefault}{\mddefault}{\updefault}{\color[rgb]{0,0,0}$a_2$}%
}}}}
\put(4276,-1411){\makebox(0,0)[lb]{\smash{{\SetFigFont{12}{14.4}{\familydefault}{\mddefault}{\updefault}{\color[rgb]{0,0,0}$a_4$}%
}}}}
\put(1351,-1411){\makebox(0,0)[lb]{\smash{{\SetFigFont{12}{14.4}{\familydefault}{\mddefault}{\updefault}{\color[rgb]{0,0,0}$b_2$}%
}}}}
\put(3676,-1411){\makebox(0,0)[lb]{\smash{{\SetFigFont{12}{14.4}{\familydefault}{\mddefault}{\updefault}{\color[rgb]{0,0,0}$b_4$}%
}}}}
\end{picture}%